\documentclass[11pt]{article}
\usepackage{amsmath,amssymb}
\newtheorem{theo}{Theorem}[section]

\newtheorem{lemma}[theo]{Lemma}
\newtheorem{coro}[theo]{Corollary}
\newtheorem{conj}[theo]{Conjecture}

\newcommand{\ignore}[1]{}
\def\square{\vrule height6pt width7pt depth1pt}
\def\endpf{\hfill\square\bigskip}

\begin{document}

\title{Mean Ramsey-Tur\'an numbers}
\author{Raphael Yuster
\thanks{
e-mail: raphy@research.haifa.ac.il \qquad
World Wide Web: http:$\backslash\backslash$research.haifa.ac.il$\backslash$\symbol{126}raphy}
\\ Department of Mathematics\\ University of
Haifa at Oranim\\ Tivon 36006, Israel}

\date{}

\maketitle
\setcounter{page}{1}
\begin{abstract}
A $\rho$-mean coloring of a graph is a coloring of the edges such that the average number of colors incident with each
vertex is at most $\rho$. For a graph $H$ and for $\rho \geq 1$, the {\em mean Ramsey-Tur\'an number}
$RT(n,H,\rho-mean)$  is the maximum number of edges a $\rho$-mean colored graph with $n$ vertices can have
under the condition it does not have a monochromatic copy of $H$.
It is conjectured that $RT(n,K_m,2-mean)=RT(n,K_m,2)$ where $RT(n,H,k)$ is 
the maximum number of edges a $k$ edge-colored graph with $n$ vertices can have
under the condition it does not have a monochromatic copy of $H$.
We prove the conjecture holds for $K_3$. We also prove that
$RT(n,H,\rho-mean) \leq RT(n,K_{\chi(H)},\rho-mean)+o(n^2)$.
This result is tight for graphs $H$ whose clique number equals their chromatic number.
In particular we get that if $H$ is a $3$-chromatic graph having a triangle then
$RT(n,H,2-mean) = RT(n,K_3,2-mean)+o(n^2)=RT(n,K_3,2)+o(n^2)=0.4n^2(1+o(1))$.
\end{abstract}

\section{Introduction}
All graphs considered are finite, undirected and simple. For standard graph-theoretic terminology see \cite{Bo}.
Ramsey and Tur\'an type problems are central problems in extremal graph theory.
These two topics intersect in {\em Ramsey-Tur\'an Theory} which is now a wide field of research with
many interesting results and open problems. The survey of Simonovits and S\'os \cite{SiSo} is an excellent
reference for Ramsey-Tur\'an Theory.

The {\em Ramsey number} $R(H,k)$ is the minimum integer $n$ such that in any $k$-coloring of the edges of
$K_n$ there is a monochromatic $H$. An edge coloring is called {\em $k$-local} if every vertex is incident
with at most $k$ colors. The {\em local Ramsey number} $R(H, k-loc)$ is the minimum integer $n$ such that 
in any $k$-local coloring of the edges of $K_n$ there is a monochromatic $H$.
An edge coloring is called {\em $\rho$-mean} if the average number of colors incident with each every vertex
is at most $\rho$. The {\em mean Ramsey number} $R(H, \rho-mean)$ is the minimum integer $n$ such that 
in any $\rho$-mean coloring of the edges of $K_n$ there is a monochromatic $H$.
Clearly, $R(H,k) \leq R(H, k-loc) \leq R(H, k-mean)$.
The relationship between these three parameters has been studied by various researchers. See, e.g.,
\cite{BoKoSc,CaTu,GyLeScTu,Sc}.
In particular, Gy\'arf\'as et. al. \cite{GyLeScTu} proved that $R(K_m,2)=R(K_m,2-loc)$.
Caro and Tuza proved that $R(K_m,2-loc)=R(K_m,2-mean)$ and
Schelp \cite{Sc} proved that $R(K_m, k-loc)=R(K_m,k-mean)$.

The {\em Ramsey-Tur\'an number} $RT(n,H,k)$ is the maximum number of edges a $k$-colored graph
with $n$ vertices can have under the condition it does not have a monochromatic copy of $H$.
We analogously define the local and mean Ramsey-Tur\'an numbers, denoted $RT(n,H,k-loc)$ and
$RT(n,H,\rho-mean)$ respectively, to be the maximum number of edges a $k$-local (resp. $\rho$-mean) colored graph
with $n$ vertices can have under the condition it does not have a monochromatic copy of $H$.
Clearly, $RT(n,H,k) \leq RT(n,H,k-loc) \leq RT(n,H,k-mean)$.

The relationship between $RT(n,H,k)$, Ramsey numbers and Tur\'an numbers is well-known.
The Tur\'an graph $T(n,k)$ is the complete $k$-partite graph with $n$ vertices whose vertex classes are as
equal as possible. Let $t(n,k)$ be the number of edges of $T(n,k)$.
Burr, Erd\H{o}s and Lov\'asz \cite{BuErLo} introduced the Ramsey function $r(H,k)$ which is the smallest integer $r$
for which there exists a complete $r$-partite graph having the property that any $k$ edge-coloring of it has a
monochromatic $H$.
For example, $r(K_m,k)=R(K_m,k)$ and $r(C_5,2)=5$.
Clearly, $RT(n,K_m,k)=t(n,R(K_m,k)-1)$. As shown in Theorem 13 in \cite{SiSo},
it follows from the Erd\H{o}s-Stone Theorem \cite{ErSt} that
$$
RT(n,H,k)=\left(1-\frac{1}{r(H,k)-1}\right){n \choose 2}+o(n^2).
$$
Clearly, a similar relationship holds between $RT(n,H,k-loc)$ and the analogous Ramsey function $r(H,k-loc)$.
However, no such relationship is known for $RT(n,H,k-mean)$. We conjecture that such a relationship holds.
\begin{conj}
\label{conj1}
$$
RT(n,H,k-mean)=\left(1-\frac{1}{r(H,k-mean)-1}\right){n \choose 2}+o(n^2).
$$
\end{conj}
Combining this with the fact that $R(K_m,2)=R(K_m,2-loc)=R(K_m,2-mean)$ we have the following stronger
conjecture for complete graphs and $k=2$.
\begin{conj}
\label{conj2}
$$
RT(n,K_m,2-mean)=RT(n,K_m,2)=t(n,R(K_m,2)-1).
$$ 
\end{conj}
For non-integral values of $\rho$ is is not even clear what the right conjecture for $RT(n,H,\rho-mean)$ should be.

The first result of this paper shows that Conjecture \ref{conj2} holds for $K_3$.
\begin{theo}
\label{triangle}
$RT(n,K_3,2-mean) = RT(n,K_3,2)= t(n,R(K_3,2)-1) = t(n,5) = \lfloor 0.4n^2 \rfloor$.
\end{theo}

The second result of this paper asserts that $RT(n,H,\rho-mean)$ is bounded by a function of the chromatic number of
$H$. In fact, for graphs whose clique number equals their chromatic number, $RT(n,H,\rho-mean)$ is essentially
determined by the chromatic number of $H$.
\begin{theo}
\label{chrom}
For all $\rho \geq 1$ and for all graphs $H$, $RT(n,H,\rho-mean) \leq RT(n,K_{\chi(H)},\rho-mean)+o(n^2)$.
In particular, if the chromatic number of $H$ equals its clique number then 
$RT(n,H,\rho-mean) = RT(n,K_{\chi(H)},\rho-mean)+o(n^2)$.
\end{theo}
The proof of Theorem \ref{chrom} uses a colored version of Szemer\'edi's Regularity Lemma
together with several additional ideas.
Notice that the trivial case $\rho=1$ in Theorem \ref{chrom} is equivalent to the Erd\H{o}s-Stone Theorem.
Combining Theorem \ref{triangle} with Theorem \ref{chrom} we obtain:
\begin{coro}
\label{coro1}
Let $H$ be a $3$-chromatic graph. Then, $RT(n,H,2-mean)  \leq 0.4n^2(1+o(1))$.
If $H$ contains a triangle then $RT(n,H,2-mean)  = 0.4n^2(1+o(1))$. \endpf
\end{coro}
The next section contains the proof of Theorem \ref{triangle}. Section 3 contains the proof of Theorem \ref{chrom}.

\section{Proof of Theorem \ref{triangle}}

We need to prove that $RT(n,K_3,2-mean) = t(n,5)$.
Since $K_5$ has a $2$-coloring with no monochromatic triangle, so does $T(n,5)$.
Hence, $RT(n,K_3,2-mean) \geq t(n,5)$.
We will show that $RT(n,K_3,2-mean) \leq t(n,5)$.
Clearly, the result is trivially true for $n < 6$, so we assume $n \geq 6$.
Our proof proceeds by induction on $n$.
Let $G$ have $n \geq 6$ vertices and more than $t(n,5)$ edges.
Clearly we may assume that $G$ has precisely $t(n,5)+1$ edges.
Consider any given $2$-mean coloring of $G$.
If $n=6$ then $G=K_6$. Recall from the introduction that $R(K_3,2-mean)=R(K_3,2)=6$.
As a $2$-mean coloring of $K_6$ contains a monochromatic triangle this base case of
the induction holds. If $n=7$ then $G$ is $K_7^-$. Again, it is trivial to check that any $2$-mean coloring of
$K_7^-$ contains a monochromatic triangle. Similarly, if $n=8$ then $G$ is a $K_8$ missing two edges and it is
straightforward to verify that any $2$-mean coloring of such a $G$ contains a monochromatic triangle.

Assume the theorem holds for all $6 \leq n' < n$ and $n \geq 9$.
For a vertex $v$, let $c(v)$ denote the number of colors incident with $v$ and let $d(v)$ denote the degree of $v$.

If some $v$ has $c(v) \geq 2$ and $d(v) \leq 4n/5$ then $G-v$ is also $2$-mean colored
and has more than $t(n-1,5)$ edges. Hence, by the induction hypothesis, $G - v$ has a monochromatic triangle.

Otherwise, if some $v$ has $c(v)=1$ and $d(v) \leq 3n/5$ then let $w$ be a vertex with maximum $c(w)$.
Then, $G - \{v,w\}$ is also $2$-mean colored
and has more than $t(n-2,5)$ edges. Hence, by the induction hypothesis, $G - \{v,w\}$ has a
monochromatic triangle.

Otherwise, if $v$ is an isolated vertex of $G$ then let $u$ and $w$ be two distinct vertices having
maximum $c(u)+c(w)$. Then, $G-\{v,u,w\}$ is $2$-mean colored
and has more than $t(n-3,5)$ edges. Hence, by the induction hypothesis, $G - \{v,u,w\}$ has a
monochromatic triangle.

We are left with the case where $\delta(G) > 3n/5$ and whenever $c(v) \geq 2$ then also $d(v) > 4n/5$.
Let $v$ be with $c(v)=1$ (if no such $v$ exists then the graph is $2$-local colored and hence contains a
monochromatic triangle as, trivially, $RT(n,K_3,2-loc)=t(n,5)$).
We may assume that $3n/5 < d(v) \leq 4n/5$, since otherwise we would have
$\delta(G) > 4n/5$ which is impossible for a graph with $t(n,5)+1$ edges.
Consider the neighborhood of $v$, denoted $N(v)$.
Clearly, if $w \in N(v)$ then $c(w) > 1$ otherwise (because $d(w) > 3n/5$) there must
be some $w' \in N(v)$ for which $(v,w,w')$ is a monochromatic triangle and we are done.
Thus, the minimum degree of $G[N(v)]$ is greater than $d(v)-n/5$.
Since $d(v) > 3n/5$ it follows that $G[N(v)]$ has minimum degree greater than
$2|N(v)|/3$. If $|N(v)|$ is divisible by $3$ then the theorem of Corr\'adi and Hajnal \cite{CoHa}
implies that $G[N(v)]$ has a triangle factor. If $|N(v)|-1$ is divisible by $3$ then
the theorem of Hajnal and Szemer\'edi \cite{HaSz} implies that $G[N(v)]$ has a factor into
$(|N(v)|-4)/3$ triangles and one $K_4$. If $|N(v)|-2$ is divisible by $3$ then,
similarly, $G[N(v)]$ has a factor into
$(|N(v)|-8)/3$ triangles and two $K_4$ or $(|N(v)|-5)/3$ triangles and one $K_5$.
Assume that $G$ has no monochromatic triangle.
The sum of colors incident with the vertices of any non-monochromatic triangle is at least $5=3 \cdot(5/3)$.
The sum of colors incident with the vertices of any $K_4$ having no monochromatic triangle is at least
$8 > 4 \cdot (5/3)$. The sum of colors incident with the vertices of any $K_5$ having no monochromatic triangle is
at least $10 > 5 \cdot (5/3)$. Thus,
$$
2 n \geq \sum_{v \in V} c(v) \geq n+ \frac{5}{3}d(v) > n+ \frac{5}{3}\cdot \frac{3}{5}n = 2n
$$
a contradiction.
\endpf

\section{Proof of Theorem \ref{chrom}}
Before we prove Theorem \ref{chrom} we need several to establish several lemmas.
\begin{lemma}
\label{adjust}
For every $\epsilon > 0$ there exists $\alpha=\alpha(\epsilon) > 0$ such that for all $m$ sufficiently large,
if a graph has $m$ vertices and more than
$RT(m,K_s,\rho-mean)+\epsilon m^2/4$ edges and is $(\rho+\alpha)$-mean colored, then it has a monochromatic
$K_s$.
\end{lemma}
{\bf Proof:}\,
Pick $\alpha$ such that $\epsilon m^2/4 > (\alpha m+1)(m-1)$ for all sufficiently large $m$.
Given a graph $G$ with $m$ vertices and more than
$RT(m,K_s,\rho-mean)+\epsilon m^2/4$ edges, consider a $(\rho+\alpha)$-mean coloring of $G$. By picking
$\lceil \alpha m \rceil$ non-isolated vertices of $G$ and deleting all edges incident with them we obtain a spanning
subgraph of $G$ with $m$ vertices, more than $RT(m,K_s,\rho-mean)+\epsilon m^2/4-(\alpha n+1)(n-1) \ge
RT(m,K_s,\rho-mean)$ edges, and which is $\rho$-mean
colored. By definition, it has a monochromatic $K_s$.
\endpf

\begin{lemma}
\label{mult}
If $n$ is a multiple of $m$ then $RT(n,K_s,\rho-mean) \geq RT(m,K_s,\rho-mean)n^2/m^2$.
\end{lemma}
{\bf Proof:}\,
Let $G$ be a graph with $m$ vertices and $RT(m,K_s,\rho-mean)$ edges having a $\rho$-mean coloring without
a monochromatic $K_s$. Let $G'$ be obtained from $G$ by replacing each vertex $v$ with an independent set
$X_v$ of size $n/m$. For $u \neq v$, we connect a vertex from $X_u$ with a vertex from $X_v$ if and only if
$uv$ is an edge of $G$, and we color this edge with the same color of $uv$.
Clearly, $G'$ has $RT(m,K_s,\rho-mean)n^2/m^2$ edges, the corresponding coloring is also $\rho$-mean,
and there is no monochromatic $K_s$ in $G'$. As $G'$ has $n$ vertices
we have that $RT(n,K_s,\rho-mean) \geq RT(m,K_s,\rho-mean)n^2/m^2$. \endpf

As mentioned in the introduction, our main tool in proving Theorem \ref{chrom}
is a colored version of Szemer\'edi's Regularity Lemma.
We now give the necessary definitions and the statement of the lemma.

Let $G=(V,E)$ be a graph, and let $A$ and $B$ be two disjoint subsets 
of $V$. If $A$ and $B$ are non-empty, let $e(A,B)$ denote the number of edges with one endpoint in $A$ and
another endpoint in $B$ and define the {\em density of edges} between $A$ and $B$ by 
$$
d(A,B) = \frac{e(A,B)}{|A||B|}.
$$
For $\gamma>0$ the pair $(A,B)$ is called {\em $\gamma$-regular} 
if for every $X \subset A$ and $Y \subset B$ satisfying $|X|>\gamma |A|$ and $ |Y|>\gamma |B|$ we have 
$$
|d(X,Y)-d(A,B)| < \gamma.
$$
An {\em equitable partition} of a set $V$ is a partition of $V$ into 
pairwise disjoint classes $V_1,\ldots,V_m$ of almost equal size, i.e., 
$\big| |V_i|-|V_j| \big| \leq 1$ for all $i,j$. An equitable partition 
of the set of vertices $V$ of $G$ into the classes $V_1,\ldots,V_m$ is
called {\em $\gamma$-regular} if $|V_i| < \gamma |V|$ for 
every $i$ and all but at most $\gamma {m \choose 2}$ of the pairs 
$(V_i,V_j)$ are $\gamma$-regular. Szemer\'edi \cite{Sz} proved the following.
 
\begin{lemma} 
\label{szemeredi}
For every $\gamma>0$, there is an integer $M(\gamma)>0$ such that for 
every graph $G$ of order $n>M$ there is a $\gamma$-regular partition of 
the vertex set of $G$ into $m$ classes, for some $1/\gamma < m < M$.
\end{lemma}
To prove Theorem \ref{chrom} we will need a colored version of the 
Regularity Lemma. Its proof is a straightforward modification of the
proof of the original result (see, e.g., \cite{KoSi} for details).

\begin{lemma}
\label{szcolor}
For every $\gamma >0$ and integer $r$, there exists an $M(\gamma,r)$ 
such that if the edges of a graph $G$ of order $n>M$
are $r$-colored $E(G) = E_1 \cup \cdots \cup E_r$, then 
there is a partition of the
vertex set $V(G)= V_1 \cup \cdots \cup V_m$, with $1/\gamma < m < M$, which is 
$\gamma$-regular simultaneously with 
respect to all graphs $G_i=(V,E_i)$ for $1 \leq i \leq r$. 
\end{lemma}

A useful notion associated with a $\gamma$-regular partition is that of a {\em 
cluster graph}. Suppose that $G$ is a graph with a $\gamma$-regular partition
$V= V_1 \cup \cdots \cup V_m$, and $\eta>0$
is some fixed constant (to be thought of as small, but much larger than $\gamma$.) The cluster
graph $C(\eta)$ is defined on the vertex set $\{1,\ldots,m\}$ by declaring $ij$ to be an edge if
$(V_i,V_j)$ is a $\gamma$-regular pair with edge density at least $\eta$.
From the definition, one might expect that if a cluster graph contains a 
copy of a fixed clique then so does the original graph. This is indeed 
the case, as established in the following well-known lemma (see 
\cite{KoSi}), which says more generally that if the cluster graph
contains a $K_s$ then, for any fixed $t$, the original graph contains the Tur\'an graph $T(st,s)$.

\begin{lemma}
\label{key}
For every $\eta>0$ and positive integers $s,t$ there exist a positive
$\gamma=\gamma(\eta,s,t)$ and a positive integer $n_0=n_0(\eta,s,t)$ with the following property.
Suppose that $G$ is a graph of order $n>n_0$ with a $\gamma$-regular partition 
$V = V_1 \cup \cdots \cup V_m$. Let $C(\eta)$ be the cluster graph of 
the partition. If $C(\eta)$ contains a $K_s$ then $G$ contains a $T(st,s)$.
\end{lemma}

\noindent
{\bf Proof of Theorem \ref{chrom}:}\,
Fix an $s$-chromatic graph $H$ and fix a real $\rho \geq 1$.
We may assume $s \geq 3$ as the theorem is trivially true (and meaningless) for bipartite graphs.
Let $\epsilon > 0$. We prove that there exists $N=N(H,\rho,\epsilon)$ such that
for all $n > N$, if $G$ is a graph with $n$ vertices and more than $RT(n, K_s,\rho-mean)+\epsilon n^2$ edges then any
$\rho$-mean coloring of $G$ contains a monochromatic copy of $H$.

We shall use the following parameters.
Let $t$ be the smallest integer for which $T(st,s)$ contains $H$.
Let $r=\lceil 18\rho^2/\epsilon^2 \rceil$.
In the proof we shall choose $\eta$ to be sufficiently small as a function of $\epsilon$ alone.
Let $\alpha=\alpha(\epsilon)$ be as in lemma \ref{adjust}.
Let $\gamma$ be chosen such that (i) $\gamma < \eta/r$, (ii) $\rho/(1-\gamma r) < \rho+\alpha$,
(iii) $1/\gamma$ is larger than the minimal $m$ for which Lemma \ref{adjust} holds.
(iv) $\gamma < \gamma(\eta,s,t)$ where $\gamma(\eta,s,t)$ is the function from Lemma \ref{key}.
In the proof we shall assume, whenever necessary, that $n$ is sufficiently large w.r.t. all of these constants, and hence
$N=N(H,\rho,\epsilon)$ exists. In particular, $N > n_0(\eta,s,t)$ where $n_0(\eta,s,t)$ is the
function from Lemma \ref{key} and also $N > M(\gamma,r)$ where $M(\gamma,r)$ is the function from
Lemma \ref{szcolor}.

Let $G=(V,E)$ be a graph with $n$ vertices and with $|E| > RT(n, K_s,\rho-mean)+\epsilon n^2$.
Notice that since $s \geq 3$ and since $RT(n, K_s,\rho-mean) \geq RT(n,K_3,1)=t(n,2)=\lfloor n^2/4 \rfloor$
we have that $n^2/2 > |E| > n^2/4$.
Fix a $\rho$-mean coloring of $G$. Assume the colors are $\{1,\ldots,q\}$ for some $q$ and let $c_i$ denote the
number of edges colored with $i$.
Without loss of generality we assume that $c_i \ge c_{i+1}$. We first show that the first $r$
colors already satisfy $c_1+c_2 + \cdots + c_r \ge |E| -\epsilon n^2/2$. Indeed, assume otherwise.
Since, trivially, $c_{r+1} \le |E|/r$, let us partition the colors $\{r+1,\ldots,q\}$ into parts such that for each part
(except, perhaps, the last part) the total number of edges colored with a color belonging to the part
is between $|E|/r$ and $2|E|/r$. The number of edges colored by a color from the last part is at most $2|E|/r$.
The number of parts is, therefore, at least
$$
\frac{\frac{\epsilon}{2}n^2}{\frac{2|E|}{r}} > \frac{\epsilon}{2}r.
$$
Since any set of $z$ edges is incident with at least $\sqrt{2z}$ vertices
we have that the total number of vertices incident with colors $r+1$ and higher is at least
$$
\left(\frac{\epsilon}{2}r -1\right)\sqrt{2|E|/r} > \frac{\epsilon}{3}r\frac{n}{\sqrt{2r}} =
\frac{\epsilon \sqrt{r}}{\sqrt{18}}n > \rho n,
$$
a contradiction to the fact that $G$ is $\rho$-mean colored.

Let $E_i$ be the set of edges colored $i$, let $G_i=(V,E_i)$, let $E'=E_1 \cup \cdots \cup E_r$ and
let $G'=(V,E')$. By the argument above, $|E'| > RT(n, K_s,\rho-mean)+\epsilon n^2/2$ and $G'$ is
$\rho$-mean colored. It suffices to show that $G'$ has a copy of $H$.

We apply Lemma \ref{szcolor} to $G'$ and obtain a partition of $V$ into $m$ classes
$V_1 \cup \cdots \cup V_m$ where $1/\gamma < m < M$ 
which is  $\gamma$-regular simultaneously with  respect to all graphs $G_i=(V,E_i)$ for $1 \le i \le r$.
Consider the cluster graph $C(\eta)$. By choosing $\eta$ sufficiently small as a function of $\epsilon$ we are
guaranteed that $C(\eta)$ has at least
$RT(m, K_s,\rho-mean)+\epsilon m^2/4$ edges. To see this, notice that if $C(\eta)$ had
less edges then, by Lemma \ref{mult}, by the definition of $\gamma$-regularity and by the definition of $C(\eta)$,
the number of edges of $G'$ would have been at most
$$
(RT(m, K_s,\rho-mean)+\frac{\epsilon}{4}m^2)\frac{n^2}{m^2} +
\eta \frac{n^2}{m^2}{m \choose 2}+\gamma{m \choose 2}
\frac{n^2}{m^2}+{{n/m} \choose 2}m
$$
$$
< RT(m, K_s,\rho-mean)\frac{n^2}{m^2}+\frac{\epsilon}{2}n^2
\leq RT(n, K_s,\rho-mean)+\frac{\epsilon}{2}n^2
$$
contradicting the cardinality of $|E'|$. In the last inequality we assume each color class has size $n/m$ precisely.
This may clearly be assumed since floors and ceilings may be dropped due to the asymptotic nature of our result.

We define a coloring of the edges of $C(\eta)$ as follows. The edge $ij$ is colored by the color whose frequency in
$E'(V_i,V_j)$ is maximal. Notice that this frequency is at least $(n^2/m^2)\eta/r$.
Let $\rho^*$ be the average number of colors incident with each vertex in this coloring of $C(\eta)$.
We will show that $\rho^* \leq \rho+\alpha$.
For $i=1,\ldots,m$ let $c(j)$ denote the number of colors incident with vertex $j$ in our coloring of $C(\eta)$.
Clearly, $c(1)+\cdots+c(m)=\rho^*m$. For $v \in V$, let $c(v)$ denote the number of colors
incident with vertex $v$ in the coloring of $G'$. Clearly, $\sum_{v \in V}c(v) \leq \rho n$.
We will show that almost all vertices $v \in V_j$ have $c(v) \geq c(j)$. Assume that color $i$ appears in vertex $j$
of $C(\eta)$. Let $V_{j,i} \subset V_j$ be the set of vertices of $V_j$ incident with color $i$ in $G'$.
We claim that $|V_j - V_{j,i}| < \gamma n/m$. Indeed, if this was not the case then by letting $Y=V_j - V_{j,i}$
and letting $X=V_{j'}$ where $j'$ is any class for which $jj'$ is colored $i$ we have that $d(X,Y)=0$ with respect to
color $i$, while $d(V_j,V_{j'}) \geq \eta/r$ with respect to color $i$. Since $\eta/r > \gamma$ this contradicts the
$\gamma$-regularity of the pair $(V_j,V_{j'})$ with respect to color $i$.
Now, let $W_j=\{v \in V_j ~:~ c(v) \geq c(j)\}$. We have therefore shown that $|W_j| \geq |V_j| - \gamma r n/m$.
Hence,
$$
\rho n \ge \sum_{v \in V}c(v) \ge \sum_{j=1}^m \sum_{v \in W_j} c(v) \ge \sum_{j=1}^m c(j)\frac{n}{m}(1-\gamma r)
=\rho^*n (1-\gamma r).
$$
It follows that
$$
\rho^* \leq \frac{\rho}{1-\gamma r} \le \rho+\alpha.
$$
We may now apply Lemma \ref{adjust} to $C(\eta)$ and obtain that $C(\eta)$ has a monochromatic $K_s$,
say with color $j$. By Lemma \ref{key} (applied to the spanning subgraph of $C(\eta)$ induced
by the edges colored $j$) this implies that $G_j=(V,E_j)$ contains a copy of $T(st,s)$.
In particular, there is a monochromatic copy of $H$ in $G$.
We have therefore proved that $RT(n,H,\rho-mean) \leq RT(n,K_s,\rho-mean)+\epsilon n^2$.
Now, if $H$ contains a $K_s$ then we also trivially have $RT(n,H,\rho-mean) \geq RT(n,K_s,\rho-mean)$.
This completes the proof of Theorem \ref{chrom}.
\endpf

\section{Acknowledgment}
The author thanks Y. Caro for useful discussions.

\end{document}